\documentclass[a4paper,11pt,twoside,reqno]{amsart}

\usepackage[utf8]{inputenc}
\usepackage[plainpages=false,pdfpagelabels=true]{hyperref}
\usepackage{amssymb,amsthm}
\usepackage[margin=1in]{geometry}
\usepackage{slashed}

\newtheorem{Satz}{Theorem}[section]

\newtheorem{Lem}[Satz]{Lemma}

\newcommand{\vol}{{\operatorname{Vol}}}
\newcommand{\tr}{{\operatorname{tr}}}
\theoremstyle{definition}

\newtheorem{Bem}[Satz]{Remark}
\newcommand{\dv}{\text{ }dV}

\allowdisplaybreaks[1]

\renewcommand{\epsilon}{\varepsilon}

\newcommand{\R}{\ensuremath{\mathbb{R}}}

\numberwithin{equation}{section}

\newcommand{\D}{\slashed{D}}
\newcommand{\p}{\slashed{\partial}}

\title[A vanishing result for the supersymmetric nonlinear sigma model]{A vanishing result for the supersymmetric nonlinear sigma model in higher dimensions}
\author{Volker Branding}
\date{\today}
\address{University of Vienna, Faculty of Mathematics\\
Oskar-Morgenstern-Platz 1, 1090 Vienna, Austria\\}
\email{volker.branding@univie.ac.at}
\subjclass[2010]{58E20; 53C27; 35J61}
\keywords{supersymmetric nonlinear sigma model; complete Riemannian manifold; Liouville theorem}
\begin{document}

\begin{abstract}
We prove a vanishing result for critical points of the supersymmetric nonlinear sigma 
model on complete non-compact Riemannian manifolds of positive Ricci curvature that admit an Euclidean type Sobolev inequality,
assuming that the dimension of the domain is bigger than two and that a certain energy is sufficiently small.
\end{abstract} 

\maketitle

\section{Introduction and results}
The supersymmetric nonlinear \(\sigma\)-model is an important model in modern quantum field theory.
The precise form of its action functional is dictated by various symmetries 
as for example the invariance under diffeomorphisms on the domain and supersymmetry transformations.
In the physics literature it is usually formulated in the language of supergeometry.
On the other hand, when neglecting the invariance under supersymmetry transformations, 
this model can also be studied using well-established tools from the geometric calculus
of variation and there is a growing interest in this mathematical version.

This mathematical study was initiated in \cite{MR2262709}, where the notion of \emph{Dirac-harmonic maps} was introduced.
The action functional for Dirac-harmonic maps consists of the well-studied Dirichlet action for a map 
between two Riemannian manifolds and the Dirac action for a \emph{vector spinor} that is defined along the map.
The critical points of this functional naturally couple the harmonic map equation with spinor fields.

However, Dirac-harmonic maps only capture part of the critical points of the supersymmetric nonlinear sigma 
model. Its full action contains additional terms:
Considering also a curvature term that is quartic in the spinors one is led to 
\emph{Dirac-harmonic maps with curvature term}, see \cite{MR3333092,MR2370260} and also \cite{MR3735550}.
Dirac-harmonic maps to target spaces with torsion are analyzed in \cite{MR3493217}.
Taking into account an additional two-form in the action functional the resulting equations were studied in \cite{MR3305429}.

Recently, there has been an intensive study of the supersymmetric nonlinear \(\sigma\)-model coupled to 
a gravitino \cite{MR3772035} which is the superpartner of the metric on the domain.

For a recent overview on the mathematical analysis of the supersymmetric nonlinear \(\sigma\)-model see the survey \cite{1710.01519}.

The action functional of the supersymmetric sigma model is special if the domain manifold is two-dimensional
since in this case the action is invariant under conformal transformations and
its critical points share special properties such as the removal of isolated singularities.

For the simplest mathematical version of the supersymmetric nonlinear \(\sigma\)-model, that is 
Dirac-harmonic maps, several existence result could be achieved.
Making use of the Atiyah-Singer index theorem uncoupled solutions of the equations for Dirac-harmonic maps
could be obtained in \cite{MR3070562}. Employing the heat-flow method an existence result for the domain
being a compact surface with boundary could be established in \cite{MR3724759}.

In this article we are aiming in the opposite direction, namely we want to give
a criterion ensuring that solutions of the critical points 
of the supersymmetric nonlinear \(\sigma\)-model must be trivial, meaning that
the map part maps to a point and the vector spinor vanishes identically.

Let us now describe the mathematical structure that we are using in more detail.
In the following we will assume that \((M,g)\) is a complete non-compact Riemannian spin manifold of dimension \(n=\dim M\geq 3\)
and \((N,h)\) a second complete Riemannian manifold of bounded geometry. Whenever we will make use of indices
we use Latin letters for indices related to \(M\) and Greek letters for indices related to \(N\).
The spin assumption guarantees the existence of the (complex) spinor bundle \(\Sigma M\) and sections
in this bundle are called \emph{spinors}. The spinor bundle is a vector bundle over the manifold \(M\) that is
equipped with a connection \(\nabla^{\Sigma M}\) and a hermitian scalar product.
On the spinor bundle there exists the algebraic
operation of Clifford multiplying a spinor with a tangent vector, which is skew-symmetric
\begin{align*}
\langle X\cdot\psi,\xi\rangle_{\Sigma M}=-\langle\psi,X\cdot\xi\rangle_{\Sigma M}
\end{align*}
for all \(X\in TM\) and \(\psi,\xi\in\Gamma(\Sigma M)\). Moreover, the Clifford relations
\begin{align*}
X\cdot Y\cdot\psi+Y\cdot X\cdot\psi=-2g(X,Y)\psi
\end{align*}
hold for all \(X,Y\in TM\) and \(\psi\in\Gamma(\Sigma M)\).

Let \(\phi\colon M\to N\) be a map and let \(\phi^\ast TN\) be the pull-back of the tangent bundle from \(N\).
We consider the twisted bundle \(\Sigma M\otimes\phi^\ast TN\), on this bundle
we obtain a connection induced from \(\Sigma M\) and \(\phi^\ast TN\), which will be
denoted by \(\tilde{\nabla}\). Sections in \(\Sigma M\otimes\phi^\ast TN\) are called \emph{vector spinors}.
On \(\Sigma M\otimes\phi^\ast TN\) we have a scalar product induced from \(\Sigma M\) and \(\phi^\ast TN\),
we will denote its real part by \(\langle\cdot,\cdot\rangle\).
The twisted Dirac operator acting on vector spinors is defined by
\[
\D:=\sum_{i=1}^ne_i\cdot\tilde{\nabla}_{e_i},
\]
where \(e_i,i=1,\ldots,\dim M\) is an orthonormal basis of \(TM\).
Note that the operator \(\D\) is elliptic. Moreover, we assume that the connection on \(\phi^\ast TN\) is metric, thus \(\D\)
is also self-adjoint with respect to the \(L^2\)-norm if \(M\) is compact.

In this article we want to focus on the system of partial differential equations that arises as
critical points of the full supersymmetric nonlinear sigma model without referring to its
concrete structure, see \cite{MR3558358} for a previous analysis.

The critical points satisfy a coupled system of the following form 
\begin{align}
\label{phi-mfd}
\tau(\phi)=&A(\phi)(d\phi,d\phi)+B(\phi)(d\phi,\psi,\psi)+C(\phi)(\psi,\psi,\psi,\psi), \\
\label{psi-mfd}
\D\psi=&E(\phi)(d\phi)\psi+F(\phi)(\psi,\psi)\psi.
\end{align}
Here, \(\tau(\phi):=\tr_g\nabla d\phi\in\Gamma(\phi^{\ast}TN)\) denotes the tension field of the map \(\phi\) and the other terms represent the analytical structure of 
the right hand side. We will always assume that the endomorphisms \(A,B,C,E\) and \(F\) are bounded.
The structure of this system is motivated by the scaling behavior of \(\phi\) and \(\psi\).
More precisely, it is the most general system one can build when ones gives weight \(0\) to the map part and weight \(1/2\) to the vector spinors.

The system \eqref{phi-mfd}, \eqref{psi-mfd} includes the cases of \emph{Dirac-harmonic maps}, which are solutions of
\begin{align*}
\tau(\phi)&=\frac{1}{2}\sum_{i=1}^nR^N(\psi,e_i\cdot\psi)d\phi(e_i), \\
\D\psi&=0
\end{align*}
and also \emph{Dirac-harmonic maps with curvature term}, which are solutions of
\begin{align}
\nonumber\tau(\phi)=&\frac{1}{2}\sum_{i=1}^nR^N(\psi,e_i\cdot\psi)d\phi(e_i)-\frac{1}{12}\langle(\nabla R^N)^\sharp(\psi,\psi)\psi,\psi\rangle, \\
\D\psi=&\frac{1}{3}R^N(\psi,\psi)\psi.
\end{align}
Here, \(R^N\) denotes the curvature tensor on \(N\).

Note that if we consider \emph{magnetic Dirac-harmonic maps} \cite{MR3305429} then \(A\neq 0\) and
in the case of \emph{Dirac-harmonic maps with torsion} \cite{MR3493217} we have \(E\neq 0\).

Recently, an existence result for a system of the form \eqref{phi-mfd}, \eqref{psi-mfd} could be obtained in the case that
the domain manifold is a Lorentzian spacetime that expands sufficiently fast \cite{MR3830277}.

In this article we will give a vanishing result for finite energy solutions of the system 
\eqref{phi-mfd}, \eqref{psi-mfd} on higher-dimensional complete Riemannian manifolds.

In the proof of our main result we will apply an \emph{Euclidean type Sobolev inequality}
of the following form
\begin{align}
\label{sobolev-inequality}
(\int_M|u|^{2n/(n-2)}\dv)^\frac{n-2}{n}\leq C_2\int_M|\nabla u|^2\dv
\end{align}
for all \(u\in W^{1,2}(M)\) with compact support,
where \(C_2\) is a positive constant that depends on the geometry of \(M\).
Such an inequality holds in \(\R^n\) and is well-known as \emph{Gagliardo-Nirenberg inequality}.
However, if one considers a non-compact complete Riemannian manifold of infinite volume
one has to make additional assumptions to ensure that an equality of the form \eqref{sobolev-inequality} 
holds.

One way of formulating these assumptions is to introduce the notion of a positive Green's function
on complete non-compact Riemannian manifolds. Let \(\Omega\) be a bounded open subset of \(M\) 
such that \(x\in\Omega\) and \(G\) be the solution of
\begin{align*}
\begin{cases}
\Delta_g=\delta_x & \text{ in } \Omega, \\ 
G=0 & \text{ on } \partial\Omega.
\end{cases}
\end{align*}
We set \(G_x^\Omega(y)=G(y)\) when \(y\in\Omega\) and \(G^\Omega_x(y)=0\) otherwise.
In addition, we set \(G_x(y)=\sup_{\Omega~~s.t.~~ x\in\Omega}G_x^\Omega(y),y\in M\).
Then we know that \cite[Theorem 3.28]{MR1481970},
\begin{enumerate}
 \item either \(G_x(y)=+\infty,~~\forall y\in M\),
 \item or \(G_x(y)<\infty,\forall y\in M\setminus\{x\}\), where \(G_x\) is the positive Green's function
 of pole \(x\).
\end{enumerate}

In the first case the manifold \(M\) is called \emph{parabolic}, 
in the second case it is called \emph{non-parabolic}.

At this point we may state \cite[Theorem 3.29]{MR1481970}.
\begin{Satz}
Let \((M,g)\) be a complete Riemannian manifold of dimension \(n\) 
that has infinite volume. Then the following two statements are equivalent
\begin{enumerate}
 \item The Euclidean type Sobolev inequality \eqref{sobolev-inequality} holds.
 \item The manifold \((M,g)\) is non-parabolic and there exists \(K>0\) such that for any \(x\in M\)
 and any \(t>0\)
 \begin{align*}
 \vol_g\big(\{y\in M ~~s.t.~~ G_x(y)>t\}\big)\leq Kt^{-n/(n-2)},
 \end{align*}
where \(G_x\) is the positive minimal Green's function of pole \(x\).
\end{enumerate}
\end{Satz}

Another way of ensuring that \eqref{sobolev-inequality} holds is the following:
If \((M,g)\) is a complete, non-compact Riemannian manifold of dimension \(n\)
with nonnegative Ricci curvature, and if for some point \(x\in M\)
\begin{align*}
\lim_{R\to\infty}\frac{\vol_g(B_R(x))}{\omega_nR^n}>0
\end{align*}
holds, then \eqref{sobolev-inequality} holds true, see \cite{MR1007511}.
Here, \(\omega_n\) denotes the volume of the unit ball in \(\R^n\).

For further geometric conditions ensuring that \eqref{sobolev-inequality} holds
we refer to \cite[Section 3.7]{MR1481970}.

In this article we will give the following result:

\begin{Satz}
\label{theorem-liouville}
Let \((M,g)\) be a complete and non-compact Riemannian manifold of dimension \(\dim M=n>2\)
with positive Ricci curvature that admits an Euclidean type Sobolev inequality of the form \eqref{sobolev-inequality}. 
Moreover, let \(N\) be a Riemannian manifold of bounded geometry.
Assume that \((\phi,\psi)\) is a smooth solution of the system \eqref{phi-mfd}, \eqref{psi-mfd}.
If
\begin{align}
\label{assumption-smallness}
\int_M(|d\phi|^{n}+|\psi|^{2n})\dv<\epsilon
\end{align}
with \(\epsilon>0\) sufficiently small, then \(\phi\) must be trivial
and \(\psi\) vanishes identically.
\end{Satz}

Note that due to a classical result of Yau \cite{MR0417452}
a complete non-compact Riemannian manifold with nonnegative Ricci curvature has infinite volume.

\begin{Bem}
In string theory one is interested in supersymmetric sigma models with 
a two-dimensional domain since one is looking for a ``generalized'' surface
of minimal area in a given target space.
However, supersymmetric sigma models also exist for domain manifolds
with dimension bigger than two, see for example \cite[Chapter 3]{MR1701600}.
Although we are neglecting supersymmetry in our analysis our 
main result Theorem \ref{theorem-liouville} may serve as a guideline 
to derive vanishing results for the supersymmetric models from physics.
\end{Bem}

\begin{Bem}
The small number \(\epsilon\) in condition \eqref{assumption-smallness} in Theorem \ref{theorem-liouville}
can also be made explicit, it has to satisfy the inequality
\begin{align*}
\nonumber
\epsilon\leq\max\{&\bigg(\frac{2n-2}{n^2C_s}\big[n|F|^2_{L^\infty}+\frac{1+n}{1+n/2}(3|C|^2_{L^\infty}+\frac{3}{2}|B|^2_{L^\infty})
+\frac{n^2}{2+n}(|R^N|_{L^\infty}+|E|^2_{L^\infty})\big]^{-1}\bigg)^\frac{n}{2}, \\
\nonumber&\bigg(\frac{2n-4}{n^2C_s}\big[|R^N|_{L^\infty}+\frac{1+n}{2}(3|A|^2_{L^\infty}+\frac{3}{2}|B|^2_{L^\infty})
+\frac{-2+n}{2+n}\frac{1+n}{2}(3|C|^2_{L^\infty}+\frac{3}{2}|B|^2_{L^\infty}) \\
&+\frac{n}{1+n/2}(|R^N|_{L^\infty}+n|E|^2_{L^\infty})]^{-1}\bigg)^\frac{n}{2}
\}.
\end{align*}
Moreover, we want to mention that demanding the smallness of the \(L^n\)-norms of \(d\phi\) and \(|\psi|^2\)
is a natural condition since these norms scale in the correct way.
\end{Bem}

\begin{Bem}
If one considers the system \eqref{phi-mfd}, \eqref{psi-mfd} 
for the domain being a closed surface several vanishing results similar to Theorem \ref{theorem-liouville}
have been obtained. More precisely, suppose that \((\phi,\psi)\) is a smooth solution of the system \eqref{phi-mfd}, \eqref{psi-mfd}
with small energy, that is
\begin{align*}
\int_M(|d\phi|^2+|\psi|^4)\dv<\epsilon.
\end{align*}
\begin{enumerate}
 \item If \((\phi,\psi)\) is a Dirac-harmonic map (\(A=C=E=F=0\)) then \(\phi\) is constant and \(\psi\)
 is a standard harmonic spinor \cite[Proposition 4.2]{MR2262709}.
 \item If \((\phi,\psi)\) is a Dirac-harmonic map with curvature term (\(A=E=0\)) and there do not exist harmonic spinors on \(M\)
 then \(\phi\) is constant and \(\psi\) vanishes identically, see \cite[Lemma 4.9]{MR3333092} and \cite[Corollary 3.7]{MR3830780}.
\end{enumerate}
\end{Bem}

This article is organized as follows: In section 2 we present the proof of the main result
and in the last section we make a short comment on Dirac-harmonic maps with curvature term
that are also critical with respect to the domain metric.

\section{Proof of theorem \ref{theorem-liouville}}
Our method of proof is inspired from a global pinching lemma established in \cite{MR2000989}
and a recent result on biharmonic maps between complete Riemannian manifolds \cite{MR3834926}.
Throughout the proof we will employ the usual summation convention, that is we will sum over repeated indices.
We will make use of a cutoff function  \(0\leq\eta\leq 1\) on \(M\) that satisfies
\begin{align*}
\eta(x)=1\textrm{ for } x\in B_R(x_0),\qquad \eta(x)=0\textrm{ for } x\in B_{2R}(x_0),\qquad |\nabla\eta|\leq\frac{C}{R}\textrm{ for } x\in M,
\end{align*}
where \(B_R(x_0)\) denotes the geodesic ball around \(x_0\) with radius \(R\).

Moreover, we employ the Weitzenböck formula for the twisted Dirac-operator \(\D\), that is
\begin{equation}
\label{weitzenboeck}
\D^2\psi=-\tilde{\Delta}\psi+\frac{1}{4}S^M\psi +\frac{1}{2}e_i\cdot e_j\cdot R^N(d\phi(e_i),d\phi(e_j))\psi.
\end{equation}
Here, \(\tilde{\Delta}\) denotes the connection Laplacian on \(\Sigma M\otimes\phi^\ast TN\), \(S^M\) denotes the scalar curvature on \(M\)
and \(R^N\) is the curvature tensor on \(N\).
This formula can be deduced from the general Weitzenböck formula for twisted Dirac operators,
see \cite[Theorem II.8.17]{MR1031992}.

\begin{Lem}
\label{lem-psi-lq}
Suppose that \(\psi\in\Gamma(\Sigma M\otimes\phi^\ast TN)\) is a smooth solution of \eqref{psi-mfd} and fix \(q>2\).
Then the following inequality holds
\begin{align}
\label{ineq-psi-lq}
\frac{C}{R^2}\int_M|\psi|^q\dv\geq&(\frac{1}{2}-\delta_1-\delta_2)\int_M\eta^2|\tilde{\nabla}\psi|^2|\psi|^{q-2}\dv
+(q-2)\int_M\eta^2|\langle\tilde{\nabla}\psi,\psi\rangle|^2|\psi|^{q-4}\dv \\
\nonumber&+\frac{1}{4}\int_M\eta^2S^M|\psi|^q\dv-C_1\int_M\eta^2|\psi|^q|d\phi|^2\dv 
-C_2\int_M\eta^2|\psi|^{q+4}\dv,
\end{align}
where \(\delta_1,\delta_2\) are two positive numbers that may be chosen arbitrarily and
\begin{align*}
C_1:=n|R^N|_{L^\infty}+\frac{n|E|^2_{L^\infty}}{4\delta_1},\qquad C_2:=\frac{n|F|^2_{L^\infty}}{4\delta_2}.
\end{align*}
\end{Lem}

\begin{proof}
Combining \eqref{psi-mfd} and \eqref{weitzenboeck} we obtain
\begin{align*}
\tilde{\Delta}\psi=&-E(\phi)(d\phi)\D\psi-F(\phi)(\psi,\psi)\D\psi-\nabla(E(\phi)(d\phi))\cdot\psi-\nabla (F(\phi)(\psi,\psi))\cdot\psi \\
&+\frac{1}{4}S^M\psi +\frac{1}{2}e_i\cdot e_j\cdot R^N(d\phi(e_i),d\phi(e_j))\psi.
\end{align*}
Testing this equation with \(\eta^2|\psi|^{q-2}\psi\) and integrating over \(M\) we find
\begin{align*}
\int_M\eta^2\langle\tilde{\Delta}\psi,\psi\rangle|\psi|^{q-2}\dv=&\frac{1}{4}\int_M\eta^2S^M|\psi|^q\dv \\
&+\frac{1}{2}\int_M\eta^2|\psi|^{q-2}\langle e_i\cdot e_j\cdot R^N(d\phi(e_i),d\phi(e_j))\psi ,\psi\rangle\dv \\
&-\int_M\eta^2|\psi|^{q-2}\langle E(\phi)(d\phi)\D\psi,\psi\rangle\dv \\
&-\int_M\eta^2|\psi|^{q-2}\langle F(\phi)(\psi,\psi)\D\psi,\psi\rangle\dv.
\end{align*}
Note that the terms involving derivatives of \(E(\phi)\) and \(F(\phi)\) vanish due to the skew-symmetry of 
the Clifford multiplication.
Moreover, using integration by parts and the properties of the cutoff function \(\eta\), we find
\begin{align*}
\int_M\eta^2\langle\tilde{\Delta}\psi,\psi\rangle|\psi|^{q-2}\dv
=&-\int_M\eta^2|\tilde\nabla\psi|^2|\psi|^{q-2}\dv
-(q-2)\int_M\eta^2|\langle\tilde\nabla\psi,\psi\rangle|^2|\psi|^{q-4}\dv\\
&-2\int_M\langle\tilde{\nabla}\psi,\psi\rangle|\psi|^{q-2}\eta\nabla\eta\dv \\
\leq &\frac{C}{R^2}\int_M|\psi|^q\dv-\frac{1}{2}\int_M\eta^2|\tilde{\nabla}\psi|^2|\psi|^{q-2}\dv \\
&-(q-2)\int_M\eta^2|\langle\tilde{\nabla}\psi,\psi\rangle|^2|\psi|^{q-4}\dv.
\end{align*}
Combining both equations we get
\begin{align*}
\frac{C}{R^2}\int_M|\psi|^q\dv\geq&\frac{1}{2}\int_M\eta^2|\tilde\nabla\psi|^2|\psi|^{q-2}\dv
+(q-2)\int_M\eta^2|\langle\tilde\nabla\psi,\psi\rangle|^2|\psi|^{q-4}\dv \\
\nonumber &+ \frac{1}{4}\int_M\eta^2S^M|\psi|^q\dv
+\frac{1}{2}\int_M\eta^2|\psi|^{q-2}\langle e_i\cdot e_j\cdot R^N(d\phi(e_i),d\phi(e_j))\psi ,\psi\rangle\dv \\
&-\int_M\eta^2|\psi|^{q-2}\langle E(\phi)(d\phi)\D\psi,\psi\rangle\dv
-\int_M\eta^2|\psi|^{q-2}\langle F(\phi)(\psi,\psi)\D\psi,\psi\rangle\dv.
\end{align*}
The result follows by using the inequalities
\begin{align*}
-\langle E(\phi)(d\phi)\D\psi,\psi\rangle&\geq -\delta_1|\tilde{\nabla}\psi|^2-\frac{-n|E|^2_{L^\infty}}{4\delta_1}|d\phi|^2|\psi|^2, \\
-\langle F(\phi)(\psi,\psi)\D\psi,\psi\rangle&\geq -\delta_2|\tilde{\nabla}\psi|^2-\frac{n|F|^2_{L^\infty}}{4\delta_2}|\psi|^6, \\
\langle e_i\cdot e_j\cdot R^N(d\phi(e_i),d\phi(e_j))\psi ,\psi\rangle&\geq 
-n|R^N|_{L^\infty}|d\phi|^2|\psi|^2,
\end{align*}
where \(\delta_1,\delta_2>0\) may be chosen arbitrarily.
\end{proof}

As a next step we derive a similar inequality for the differential of the map \(\phi\).
\begin{Lem}
\label{lem-phi-lr}
Suppose that \(\phi\colon M\to N\) is a smooth solution of \eqref{phi-mfd} and fix \(r>2\).
Then the following inequality holds
\begin{align}
\label{ineq-phi-lr}
\frac{C}{R^2}\int_M|d\phi|^r\dv\geq&\int_M\eta^2\langle d\phi(\text{Ric}^M(e_i)),d\phi(e_i)\rangle|d\phi|^{r-2}\dv 
+\frac{1}{2}\int_M\eta^2|\nabla d\phi|^2|d\phi|^{r-2}\dv \\
\nonumber&+\frac{r-2}{2}\int_M\eta^2|\langle\nabla d\phi,d\phi\rangle|^2|d\phi|^{r-4}\dv \\
\nonumber&-C_3\int_M\eta^2|d\phi|^{r+2}\dv
-C_4\int_M\eta^2|\psi|^8|d\phi|^{r-2}\dv,
\end{align}
where 
\begin{align*}
C_3:=&|R^N|_{L^\infty}+\frac{1+r}{2}(3|A|^2_{L^\infty}+\frac{3}{2}|B|^2_{L^\infty}), \qquad
C_4:=\frac{1+r}{2}(3|C|^2_{L^\infty}+\frac{3}{2}|B|^2_{L^\infty}).
\end{align*}
\end{Lem}
\begin{proof}
Recall that for a map between two Riemannian manifolds the following Bochner formula holds (\cite[Proposition 1.34]{MR1391729})
\begin{align*}
(\Delta d\phi)(e_i)=d\phi(\text{Ric}^M(e_i))+R^N(d\phi(e_j),d\phi(e_i))d\phi(e_j)+\nabla_{e_i}\tau(\phi).
\end{align*}
Testing this equation with \(\eta^2|d\phi|^{r-2}d\phi(e_i)\) and integrating over \(M\) we find
\begin{align*}
\int_M\eta^2\langle\Delta d\phi,d\phi\rangle|d\phi|^{r-2}\dv
=&\int_M\eta^2\langle d\phi(\text{Ric}^M(e_i)),d\phi(e_i)\rangle|d\phi|^{r-2}\dv \\
&+\int_M\eta^2\langle R^N(d\phi(e_i),d\phi(e_j))d\phi(e_i),d\phi(e_j)\rangle|d\phi|^{r-2}\dv\\
&+\int_M\eta^2\langle\nabla\tau(\phi),d\phi\rangle|d\phi|^{r-2}\dv.
\end{align*}

Using integration by parts we may rewrite
\begin{align*}
\int_M\eta^2\langle\nabla\tau(\phi),d\phi\rangle|d\phi|^{r-2}\dv=&-2\int_M\eta\nabla\eta\langle\tau(\phi),d\phi\rangle|d\phi|^{r-2}\dv
-\int_M\eta^2|\tau(\phi)|^2|d\phi|^{r-2}\dv \\
&-(r-2)\int_M\eta^2\langle\tau(\phi),d\phi\rangle|d\phi|^{r-4}\langle d\phi,\nabla d\phi\rangle\dv
\end{align*}
and also
\begin{align*}
\int_M\eta^2\langle\Delta d\phi,d\phi\rangle|d\phi|^{r-2}\dv=&-2\int_M\eta\nabla\eta\langle\nabla d\phi,d\phi\rangle|d\phi|^{r-2}\dv 
-\int_M\eta^2|\nabla d\phi|^2|d\phi|^{r-2}\dv \\
&-(r-2)\int_M\eta^2|\langle\nabla d\phi,d\phi\rangle|^2|d\phi|^{r-4}\dv.
\end{align*}
This allows us to deduce the following inequality
\begin{align*}
\int_M\eta^2&\langle d\phi(\text{Ric}^M(e_i)),d\phi(e_i)\rangle|d\phi|^{r-2}\dv
+\int_M\eta^2|\nabla d\phi|^2|d\phi|^{r-2}\dv \\
&+(r-2)\int_M\eta^2|\langle\nabla d\phi,d\phi\rangle|^2|d\phi|^{r-4}\dv \\
=&2\int_M\eta\nabla\eta\langle\tau(\phi),d\phi\rangle|d\phi|^{r-2}\dv
-2\int_M\eta\nabla\eta\langle\nabla d\phi,d\phi\rangle|d\phi|^{r-2}\dv \\
&-\int_M\eta^2\langle R^N(d\phi(e_i),d\phi(e_j))d\phi(e_i),d\phi(e_j)\rangle|d\phi|^{r-2}\dv
+\int_M\eta^2|\tau(\phi)|^2|d\phi|^{r-2}\dv \\
&+(r-2)\int_M\eta^2\langle\tau(\phi),d\phi\rangle|d\phi|^{r-4}\langle d\phi,\nabla d\phi\rangle\dv \\
\leq &\frac{C}{R^2}\int_M|d\phi|^r\dv+\frac{1}{2}\int_M\eta^2|\nabla d\phi|^2|d\phi|^{r-2}\dv 	
+\frac{1+r}{2}\int_M\eta^2|\tau(\phi)|^2|d\phi|^{r-2}\dv \\
&+|R^N|_{L^\infty}\int_M\eta^2|d\phi|^{r+2}\dv
+\frac{r-2}{2}\int_M\eta^2|\langle\nabla d\phi,d\phi\rangle|^2|d\phi|^{r-4}\dv,
\end{align*}
which yields
\begin{align*}
\int_M\eta^2&\langle d\phi(\text{Ric}^M(e_i)),d\phi(e_i)\rangle|d\phi|^{r-2}\dv
+\frac{1}{2}\int_M\eta^2|\nabla d\phi|^2|d\phi|^{r-2}\dv \\
&+\frac{r-2}{2}\int_M\eta^2|\langle\nabla d\phi,d\phi\rangle|^2|d\phi|^{r-4}\dv \\
\leq &\frac{C}{R^2}\int_M|d\phi|^r\dv
+\frac{1+r}{2}\int_M\eta^2|\tau(\phi)|^2|d\phi|^{r-2}\dv 
+|R^N|_{L^\infty}\int_M\eta^2|d\phi|^{r+2}\dv.
\end{align*}

Moreover, making use of \eqref{phi-mfd} we may estimate
\begin{align*}
|\tau(\phi)|^2\leq (3|A|^2_{L^\infty}+\frac{3}{2}|B|^2_{L^\infty})|d\phi|^4+(3|C|^2_{L^\infty}+\frac{3}{2}|B|^2_{L^\infty})|\psi|^8.
\end{align*}
The result follows by combining both equations.
\end{proof}

We may combine Lemmas \ref{lem-psi-lq} and \ref{lem-phi-lr} to obtain the following 
\begin{Lem}
Suppose that the pair \((\phi,\psi)\) is a smooth solution of \eqref{phi-mfd}, \eqref{psi-mfd}.
Then the following inequality holds
\begin{align}
\label{combined-inequality}
\nonumber\frac{C}{R^2}\int_M(|\psi|^{4p}+|d\phi|^{2p})\dv\geq&\int_M\eta^2\langle d\phi(\text{Ric}^M(e_i)),d\phi(e_i)\rangle|d\phi|^{2p-2}\dv+\frac{1}{4}\int_M\eta^2S^M|\psi|^{4p}\dv \\
\nonumber&+\frac{2p-1}{2p^2}\int_M\eta^2\big|d|\psi|^{2p}\big|^2\dv
+\frac{p-1}{p^2}\int_M\eta^2\big|d|d\phi|^p\big|^2\dv \\
\nonumber&+(\frac{1}{2}-\delta_1-\delta_2)\int_M\eta^2|\tilde{\nabla}\psi|^2|\psi|^{4p-2}\dv \\
\nonumber&+\frac{1}{2}\int_M\eta^2|\nabla d\phi|^2|d\phi|^{2p-2}\dv \\
\nonumber&-(C_2+\frac{2C_4}{1+p}+\frac{pC_1}{1+p})\int_M\eta^2|\psi|^{4p+4}\dv \\
&-(C_3+C_4\frac{-1+p}{1+p}+\frac{C_1}{1+p})\int_M\eta^2|d\phi|^{2p+2}\dv,
\end{align}
where \(p>1\).
\end{Lem}
\begin{proof}
First, we add up \eqref{ineq-psi-lq} with \(q=4p\) and \eqref{ineq-phi-lr} with \(r=2p\) to get the following inequality
\begin{align*}
\frac{C}{R^2}\int_M(|\psi|^{4p}+|d\phi|^{2p})\dv\geq&\int_M\eta^2\langle d\phi(\text{Ric}^M(e_i)),d\phi(e_i)\rangle|d\phi|^{2p-2}\dv+\frac{1}{4}\int_M\eta^2S^M|\psi|^{4p}\dv \\
&+(4p-2)\int_M\eta^2|\langle\tilde\nabla\psi,\psi\rangle|^2|\psi|^{4p-4}\dv \\
&+(\frac{1}{2}-\delta_1-\delta_2)\int_M\eta^2|\tilde{\nabla}\psi|^2|\psi|^{4p-2}\dv \\
&+(p-1)\int_M\eta^2|\langle\nabla d\phi,d\phi\rangle|^2|d\phi|^{2p-4}\dv \\ 
&+\frac{1}{2}\int_M\eta^2|\nabla d\phi|^2|d\phi|^{2p-2}\dv\\
&-C_1\int_M\eta^2|\psi|^{4p}|d\phi|^2\dv-C_2\int_M\eta^2|\psi|^{4p+4}\dv \\
&-C_3\int_M\eta^2|d\phi|^{2p+2}\dv-C_4\int_M\eta^2|\psi|^8|d\phi|^{2p-2}\dv.
\end{align*}

In order to estimate the terms on the right hand side containing both \(\psi\) and \(d\phi\) we use 
the general Young inequality to obtain
\begin{align*}
|\psi|^8|d\phi|^{2p-2}\leq&\frac{2}{1+p}|\psi|^{4p+4}+\frac{-1+p}{1+p}|d\phi|^{2p+2}, \\
|d\phi|^2|\psi|^{4p}\leq&\frac{1}{1+p}|d\phi|^{2p+2}+\frac{p}{1+p}|\psi|^{4p+4}.
\end{align*}
In addition, we have
\begin{align*}
|\langle\tilde{\nabla}\psi,\psi\rangle|^2|\psi|^{4p-4}
=&\frac{1}{4}\big|d|\psi|^2\big|^2|\psi|^{4p-4}
=\frac{1}{4p^2}\big|d|\psi|^{2p}\big|^2,\\
|\langle\nabla d\phi,d\phi\rangle|^2|d\phi|^{2p-4}
=&\frac{1}{4}\big|d|d\phi|^2\big|^2|d\phi|^{2p-4}
=\frac{1}{p^2}\big|d|d\phi|^p\big|^2.
\end{align*}
Combining all the estimates then yields the claim.
\end{proof}

Making use of \eqref{sobolev-inequality} we can now give the following
estimate:

\begin{Lem}
\label{lem-sobolev}
Let \((M,g)\) be a complete non-compact Riemannian manifold of infinite volume that admits 
an Euclidean type Sobolev inequality.
Assume that \(f\) is a positive function on \(M\).
For \(\dim M=n>2\) and \(p>1\) arbitrary the following inequality holds
\begin{align}
\label{inequality-sobolev}
\int_M\eta^2f^{2p+2}\dv\leq C_s\big(\int_Mf^n\dv\big)^\frac{2}{n}\big(\frac{1}{R^2}\int_Mf^{2p}\dv
+\int_M\eta^2|df^p|^2\dv\big),
\end{align}
where the positive constant \(C_s\) depends on \(n,p\) and the geometry of \(M\).
\end{Lem}

\begin{proof}
By Hölder's inequality we find
\begin{align*}
\int_M\eta^2f^{2p+2}\dv\leq\big(\int_M(\eta f^p)^{2r}\dv\big)^\frac{1}{r}\big(\int_Mf^{\frac{2r}{r-1}}\dv\big)^{\frac{r-1}{r}}.
\end{align*}
Now, we choose \(r=\frac{n}{n-2}\) and apply \eqref{sobolev-inequality} to obtain
\begin{align*}
\big(\int_M(\eta f^p)^{\frac{2n}{n-2}}\dv\big)^\frac{n-2}{n}\leq C_2\int_M|d(\eta f^p)|^2\dv.
\end{align*}
Hence, we find
\begin{align*}
\int_M\eta^2f^{2p+2}\dv\leq C\big(\int_Mf^n\dv\big)^\frac{2}{n}\int_M|d(\eta f^p)|^2\dv.
\end{align*}
Using the properties of the cutoff function \(\eta\) we find
\begin{align*}
\int_M|d(\eta f^p)|^2\dv\leq\frac{C}{R^2}\int_Mf^{2p}\dv+2\int_M\eta^2|df^p|^2\dv,
\end{align*}
which completes the proof.
\end{proof}

At this point we are ready to give the proof of the main result.
\begin{proof}[Proof of Theorem \ref{theorem-liouville}]
We make use of the inequality \eqref{combined-inequality}, where we choose \(\delta_1=\delta_2=\frac{1}{4}\).
Moreover, we apply \eqref{inequality-sobolev} to \(f=|d\phi|,|\psi|^2\) taking into account the smallness condition \eqref{assumption-smallness}.
Choosing \(p=n/2\) we obtain the following inequality
\begin{align*}
\nonumber\frac{C}{R^2}\int_M(|\psi|^{2n}&+|d\phi|^{n})\dv \\
\geq&\int_M\eta^2\langle d\phi(\text{Ric}^M(e_i)),d\phi(e_i)\rangle|d\phi|^{n-2}\dv+\frac{1}{4}\int_M\eta^2S^M|\psi|^{2n}\dv \\
\nonumber&+\bigg(\frac{2n-2}{n^2}-\epsilon^\frac{2}{n}C_s(C_2+\frac{2C_4}{1+n/2}+\frac{nC_1}{2+n})\bigg)\int_M\eta^2\big|d|\psi|^{n}\big|^2\dv \\
\nonumber&+\bigg(\frac{2n-4}{n^2}-\epsilon^\frac{2}{n}C_s(C_3+C_4\frac{-2+n}{2+n}+\frac{C_1}{1+n/2})\bigg)\int_M\eta^2\big|d|d\phi|^\frac{n}{2}\big|^2\dv \\
\nonumber&-\epsilon^\frac{2}{n}C_s(C_2+\frac{2C_4}{1+n/2}+\frac{nC_1}{2+n})\frac{1}{R^2}\int_M|\psi|^{2n}\dv \\
&-\epsilon^\frac{2}{n}C_s(C_3+C_4\frac{-2+n}{2+n}+\frac{C_1}{1+n/2})\frac{1}{R^2}\int_M|d\phi|^{n}\dv.
\end{align*}
Choosing \(\epsilon\) small enough 
and taking the limit \(R\to\infty\) we obtain
\begin{align*}
0\geq&\int_M\langle d\phi(\text{Ric}^M(e_i)),d\phi(e_i)\rangle|d\phi|^{n-2}\dv+\frac{1}{4}\int_MS^M|\psi|^{2n}\dv.
\end{align*}
By assumption \(M\) has positive Ricci curvature and thus \(\text{Ric}^M\) becomes a positive definite 
non-degenerate bilinear form on \(TM\). Moreover, the positivity of the Ricci curvature also implies that the scalar 
curvature is positive, which completes the proof.
\end{proof}

\section{A remark on Dirac-harmonic maps with curvature term in higher dimensions}
In this section we want to make a comment on Dirac-harmonic maps with curvature term from
domain manifolds with \(\dim M\geq 3\).
The action functional for Dirac-harmonic maps with curvature term is given by
\begin{equation}
\label{energy-dhc}
S(\phi,\psi)=\frac{1}{2}\int_M(|d\phi|^2+\langle\psi,\D\psi\rangle-\frac{1}{6}\langle R^N(\psi,\psi)\psi,\psi\rangle)\dv.
\end{equation}
Here, the indices are contracted as 
\[
\langle R^N(\psi,\psi)\psi,\psi\rangle=\sum_{\alpha,\beta,\gamma,\delta}R_{\alpha\beta\gamma\delta}\langle\psi^\alpha,\psi^\gamma\rangle\langle\psi^\beta,\psi^\delta\rangle,
\]
which ensures that the functional is real-valued.
The critical points of the action functional \eqref{energy-dhc} are given by
\begin{align}
\label{euler-lagrange-phi}\tau(\phi)=&\frac{1}{2}\sum_{i=1}^nR^N(\psi,e_i\cdot\psi)d\phi(e_i)
-\frac{1}{12}\langle(\nabla R^N)^\sharp(\psi,\psi)\psi,\psi\rangle, \\
\label{euler-lagrange-psi}\D\psi=&\frac{1}{3}R^N(\psi,\psi)\psi,
\end{align}
where \(\sharp\colon\phi^\ast T^\ast N\to\phi^\ast TN\) represents the musical isomorphism.
For a derivation of the critical points 
see \cite[Section II]{MR2370260} and \cite[Proposition 2.1]{MR3333092}.

Solutions \((\phi,\psi)\) of the system \eqref{euler-lagrange-phi}, \eqref{euler-lagrange-psi}
are called \emph{Dirac-harmonic maps with curvature term}. It is obvious that Dirac-harmonic maps
with curvature term have the analytical structure of the system \eqref{phi-mfd}, \eqref{psi-mfd}.

In physics one often considers the metric on the domain of the action functional \eqref{energy-dhc} as 
a dynamic field of the theory. For this reason one also varies the action functional
with respect to the domain metric. As the corresponding Euler-Lagrange equation one
gets the vanishing of the energy-momentum tensor
\begin{align}
\label{euler-lagrance-metric}
0=S_{ij}=&2\langle d\phi(e_i),d\phi(e_j)\rangle-g_{ij}|d\phi|^2
+\frac{1}{2}\langle\psi,e_i\cdot\tilde\nabla_{e_j}\psi+e_j\cdot\tilde\nabla_{e_i}\psi\rangle 
-\frac{1}{6}g_{ij}\langle R^N(\psi,\psi)\psi,\psi\rangle.
\end{align}
For a derivation of the energy-momentum tensor see \cite[Lemma 5.3]{1605.03453}.
In the mathematics literature one calls an action functional that is critical
with respect to the domain metric \emph{stationary}.
A vanishing result for stationary Dirac-harmonic maps with curvature term 
with finite energy from complete Riemannian manifolds has been obtained in \cite[Theorem 5.4]{1605.03453}.

For a smooth Dirac-harmonic map with curvature term that is also critical with respect
to the domain metric, that is a solution of the system \eqref{euler-lagrange-phi}, \eqref{euler-lagrange-psi}, \eqref{euler-lagrance-metric},
we obtain the following vanishing result
\begin{Satz}
Let the triple \((\phi,\psi,g)\) be a smooth Dirac-harmonic map with curvature term that is critical with respect to the domain
metric. Suppose that there do not exist harmonic spinors on \(M\), that is solutions of \(\p\Psi=0\),
where \(\p\) is the standard Dirac operator and \(\Psi\in\Gamma(\Sigma M)\).

If \(\dim M\geq 3\) and the target has positive sectional curvature then
\(\phi\) is constant and \(\psi\) vanishes identically.
\end{Satz}
\begin{proof}
Taking the trace of \eqref{euler-lagrance-metric} and using \eqref{euler-lagrange-psi} we find
\begin{align*}
0=&(2-n)|d\phi|^2+\langle\psi,\D\psi\rangle-\frac{n}{6}\langle R^N(\psi,\psi)\psi,\psi\rangle \\
=&(2-n)(|d\phi|^2+\frac{1}{6}\langle R^N(\psi,\psi)\psi,\psi\rangle).
\end{align*}
At this point we need to distinguish two cases.
\begin{enumerate}
 \item First, suppose that \(\phi\colon M\to N\) is a constant map.
 Consider \(v\in\phi^\ast TN,\Psi\in\Gamma(\Sigma M)\) and set \(\psi:=\Psi\otimes v\).
 It is easy to check that this pair \((\phi,\psi)\) satisfies
 \begin{align*}
  \langle R^N(\psi,\psi)\psi,\psi\rangle=\langle R^N(v,v)v,v\rangle|\Psi|^4=0
 \end{align*}
 due to the skew symmetry of the Riemann curvature tensor regardless of any curvature assumptions on the target.
 Since \(d\phi=0\) the equations
 \eqref{euler-lagrange-phi}, \eqref{euler-lagrange-psi} reduce to \(\p\Psi=0\).
 But as we are assuming that \(M\) does not admit harmonic spinors we can deduce that \(\psi=0\).
 
 \item Now, we are considering a pair \((\phi,\psi)\) that is not of the form from above, in this case
 the term \(\langle R^N(\psi,\psi)\psi,\psi\rangle\) will be different from zero.
 A careful inspection reveals that for \(n>2\) and \(N\) having positive sectional curvature we have 
\begin{align*}
|d\phi|^2+\frac{1}{6}\langle R^N(\psi,\psi)\psi,\psi\rangle\geq 0,
\end{align*}
see \cite[Proof of Theorem 1.2]{MR2370260} for more details.

This allows us to conclude that under the assumptions of the theorem the map \(\phi\) is constant and the vector spinor \(\psi\) is trivial.
Note that we actually do not require that the map \(\phi\) satisfies \eqref{euler-lagrange-phi} here.

\end{enumerate}

\end{proof}

\begin{Bem}
It is clear that such a result does not hold if \(n=2\) since the action functional
\eqref{energy-dhc} is conformally invariant in this case.
The above result would also hold if we consider an additional two-form in the action functional as in \cite{MR3305429}
since this does not give a contribution to the energy-momentum tensor.
\end{Bem}

\par\medskip
\textbf{Acknowledgements:}
The author gratefully acknowledges the support of the Austrian Science Fund (FWF) 
through the project P30749-N35 ``Geometric variational problems from string theory''.
\bibliographystyle{plain}
\bibliography{mybib}
\end{document}